\documentclass[reqno]{amsart}
\theoremstyle{plain}
\newtheorem{thm}{\bf Theorem}[section]
\newtheorem{prop}[thm]{\bf Proposition}

\newtheorem{lem}[thm]{\bf Lemma}
\theoremstyle{remark}
\newtheorem{defn}[thm]{\bf Def{}inition}
\newtheorem{rem}[thm]{\bf Remark}
\newtheorem{exa}[thm]{\bf Example}
\numberwithin{equation}{section}
\usepackage{enumerate}

\usepackage {pstricks}
\begin{document}
\baselineskip15pt

\title[On Retro Frame Associated With Measurable Space]{On Retro Frame Associated With Measurable Space}

\author[Raj Kumar, Ashok K. Sah, Satyapriya and Sheetal]{Raj Kumar$^1$, Ashok K. Sah$^2$, Satyapriya$^3$ and Sheetal$^4$}

\address{$^{1}$ Department of Mathematics, Kirori Mal College,
University of Delhi, Delhi-110007, India.}
\email{\textcolor[rgb]{0.00,0.00,0.84}{rajkmc@gmail.com}}

\address{$^{2}$ Department of Mathematics, University of Delhi, Delhi-110007, India.}
\email{\textcolor[rgb]{0.00,0.00,0.84}{ashokmaths2010@gmail.com}}

\address{$^{3}$ Department of Mathematics, University of Delhi, Delhi-110007, India.}
\email{\textcolor[rgb]{0.00,0.00,0.84}{kmc.satyapriya@gmail.com}}

\address{$^{4}$ Department of Mathematics, SGTB Khalsa College, University of Delhi, Delhi.}
\email{\textcolor[rgb]{0.00,0.00,0.84}{sheetal@sgtbkhalsa.du.ac.in}}

\begin{abstract}\baselineskip10pt
Frames are redundant system which are useful in the reconstruction
of certain classes of spaces. Duffin and Schaeffer introduced frames for Hilbert spaces,
while addressing some deep problems in nonharmonic Fourier series. The dual of a frame (Hilbert) always
exists and can be obtained in a natural way. In this paper we introduce the notion of
$\Omega_0$-type duality of retro $(\Omega,\mu)$-frames are given.  Necessary and
sufficient conditions for the existence of the dual of retro $(\Omega,\mu)$-
frames are obtained. A special class of retro $(\Omega,\mu)$-frames which
always admit a dual frame is discussed.
\end{abstract}

\keywords{Frames, $(\Omega,\mu)$-frames, Retro $(\Omega,\mu)$-frames. 
\newline \indent $^{1}$ Corresponding author. Email: rajkmc@gmail.com}

\maketitle \thispagestyle{empty} \baselineskip13pt
\section{Introduction}

Frames for Hilbert space were first  introduced by Duffin and
Schaeffer in \cite{DS}, while working in deep problem in non-harmonic
Fourier series.  For some reason the work of Duffin and Schaeffer
was not continued until 1986, when the fundamental work of
Daubechies, Grossmann and Mayer \cite{DGM} brought this all back to life,
right at the dawn of the ``wavelet era". Recently, generalization of frames was proposed by Kaiser \cite{K} and independently by Ali Tawreque, Antoine and Gazedu \cite{AAG} who named it as a continues frames while Kaiser used the terminology generalized frames. Gabardo and Han \cite{GH} studied continues frames and use the terminology $(\Omega,\mu)$-frames.
Discrete and continuous frames arise in many applications in both pure and applied mathematics and, in particular, they play important roles in digital signal processing and scientific computations.For a nice introduction to frames an interested reader may refer to \cite{CK,PC,C,H,HH,K,KV,Y1} and references therein.
 While studying Banach frames in Banach spaces, Jain, Kaushik and Vashisht \cite{J,JK}
introduced the notion of retro Banach frames in conjugate Banach spaces. Motivated by there work
in this paper we introduce the notion of dual (or simply
$\Omega_0$-dual) of retro $(\Omega,\mu)$-frames in measurable spaces. It is observed that
even an exact retro $(\Omega,\mu)$-frame does not admit a dual retro $(\Omega,\mu)$-frame for the underlying space.
The said situation is entirely
different from the dual of frames (Hilbert) for Hilbert spaces,
where an exact frame admit a unique dual. Necessary and sufficient
conditions for the existence of retro dual $(\Omega,\mu)$-frames are obtained.
Counterexamples are also  given to defend results. A
special class of retro $(\Omega,\mu)$-frames which always admit retro dual
$(\Omega,\mu)$-frames is discussed.

\section{Preliminaries}
\textbf{Duality of frames in Hilbert spaces}: Let us
give now a   brief discussion on duality of frames in Hilbert
spaces. Suppose that $\mathcal{E} \equiv \{f_k\}$ is a frame for a
Hilbert space $\mathcal{H}$. A system $\mathcal{E}' \equiv \{g_k\}
\subset \mathcal{H}$ is called a \emph{dual frame}  of $\mathcal{E}$
if $\mathcal{E}'$ is a frame for $\mathcal{H}$ and
\begin{align*}
f=\sum_{k=1}^{\infty}\langle f, g_k \rangle f_k, \ \text{for  all}
\ f \in \mathcal{H}.
\end{align*}
In short we say that $(\mathcal{E}, \mathcal{E}')$ is a dual pair.
Let $\mathcal{E} \equiv \{f_k\}$ be a frame for
$\mathcal{H}$ with frame bounds $0 < A, \  B < \infty$ and with frame operator $S$. Then
$\mathcal{E}_{(natural)} \equiv \{S^{-1}f_k\}$ is  a frame for $\mathcal{H}$ with frame bounds $A^{-1},
B^{-1}$.The system  $\mathcal{E}_{(natural)} \equiv \{S^{-1}f_k\}$ is a dual of
$\mathcal{E} \equiv \{f_k\}$. The dual frame $\mathcal{E}_{(natural)} \equiv \{S^{-1}f_k\}$  is known as
\emph{canonical dual} (or \emph{natural dual}) of $\mathcal{E} \equiv \{f_k\}$ . Thus, every frame for a Hilbert space admit a
dual (at least one).  An interesting property of dual frames is that dual of a
frame need not be unique. In fact, there may be infinitely many
duals for a non-exact frame. A frame has a unique dual if and only
if it is exact. For more technical details about duality of frames (Hilbert), one may refer to [5, 12].\\
The following example shows that the dual of a frame (Hilbert) need not
be unique. Infact, it may have infinitely many duals.

\begin{exa}
Let $\{\chi_n\}$ be an orthogonal basis for a Hilbert space
$\mathcal{H}$. Then, $\{\chi_n\}$ is a Parseval frame for
$\mathcal{H}$. Consider a  system $\mathcal{E} \equiv \{f_k\} = \{
\chi_1, \chi_1, \chi_2, \chi_3, \chi_4,.....\}$. One can easily
verify that
\begin{align*}
\|f\|^2\leq  \|\{\langle f, f_k\rangle\} \|^2_{\ell^2} \leq
2\|f\|^2, \ \text{for \ all} \  f \in \mathcal{H}.
\end{align*}
Therefore, $\mathcal{E}$ is a frame for $\mathcal{H}$ with one of
the choice of  bounds $A=1, B=2$. Let $S$ be the frame operator of
$\mathcal{E}$. Then, the canonical dual of $\mathcal{E}$ is given by
\begin{align*}
\mathcal{E}'\equiv \{ S^{-1}(f_k)\} = \{\frac{1}{2}\chi_1,
\frac{1}{2}\chi_1, \chi_2, \chi_2,....\} .
\end{align*}
Other dual of $\mathcal{E}$ are $\mathcal{E}''\equiv \{0, \chi_1, 0,
\chi_2, .....\}$ and $\mathcal{E}'''\equiv
\{\frac{1}{3}\chi_1,\frac{2}{3}\chi_1, \chi_2, \chi_3,
\chi_4,.....\}$. Therefore, infinitely many dual of a frame
(Hilbert) can be constructed.
\end{exa}
Now we give basic definitions and notations which will be required
in this paper. Throughout this paper $\mathcal{H}$ will denote an
infinite dimensional Hilbert space over the field $\mathbb{K}$
$(\mathbb{K}=\mathbb{R}$ or $\mathbb{C})$,and $\mathcal{Z}_d$ is a Hilbert
space of scalar-valued sequences indexed by $\Omega$ which is
associated with $\mathcal{H}$. For a sequence $\{F(w)\}_{w\in \Omega}$  in
$\mathcal{H}$, $[F(w)]$ denotes the closure of the linear hull of
$\{F(w)\}_w\in \Omega$ in the norm topology of $\mathcal{H}$. As usual
$\delta_{w_0,w}$ denote the Kronecker delta which is defined as
$\delta_{w_0,w}=0,$ if $w_0 \neq w$ and  $\delta_{w_0,w}=1,$ if $w_0 = w$.
The class of all bounded linear operator from a Hilbert space
$\mathcal{H}$ into a Hilbert space $\mathcal{H}$ is denoted by
$\mathfrak{B}(\mathcal{H})$. Unless otherwise stated
all systems of the form $\{F(w)\}$ are indexed by $\Omega$.
\begin{defn}
Let $(\Omega,\mu)$ be a measure space and $\mathcal{H}$ be a Hilbert space with inner product. A vector-value mapping $F:\Omega\rightarrow \mathcal{H}$ (i.e a collection of vectors $\mathcal{F}\equiv \{F(w)\}_{w\in \Omega})$ is said to be a $(\Omega, \mu)$- frame for $\mathcal{H}$ if\\
\begin{enumerate}[(i)]
\item F is a weakly measurable function.
\item There exist constants $A$ and $B$ with $0<A\leq B$ such that
\begin{align}
A\|x\|^2\leq\int_\Omega|\langle x,F(w)\rangle|^2d\mu(w)\leq B\|x\|^2, \ \text{for  all} \ x\in \mathcal{H}
\end{align}
\end{enumerate}
\end{defn}
The positive constants $A$ and $B$, respectively, are are called lower and upper frame bounds of the $(\Omega, \mu)$-frame $\mathcal{F}\equiv \{F(w)\}_{w\in \Omega}$. They are not unique. The inequality (2.1) is called the $(\Omega, \mu)$-\textit{frame inequality}. If $A=B$, then $\{F(w)\}_{w\in \Omega}$ is called \textit{tight} and \textit{normalized tight} if $A=B=1$.The positive constants
\begin{align*}
& \mathrm{\alpha}_0= \inf \{B: B \text{ satisfy (2.1)}\}\\
&\mathrm{\beta}_0= \sup \{A: A \text{ satisfy (2.1)}\}
\end{align*}
are called \emph{optimal} or \emph{best bounds} for the   $(\Omega, \mu)$-frame $\{F(w)\}_{w\in \Omega}$. A $(\Omega, \mu)$-\textit{frame inequality} frame is said to be \textit{exact} if for arbitarary $\Omega_0\subset\Omega$, with $\mu(\Omega_0)>0$, $\{F(w)\}_{w\in \Omega\sim\Omega_0}$ ceases to be a frame for $\mathcal{H}$. If only upper inequality of (2.1) holds then $\mathcal{F}\equiv \{F(w)\}_{w\in \Omega}$ is a $(\Omega,\mu)$-\textit{Bessel family}.

 The operator $T_F:\mathcal{H}\rightarrow L^2(\Omega,\mu)$ given by
\begin{align*}
(T_Fx)(w)=\langle x,F(w)\rangle,\ w\in\Omega, \ x\in \mathcal{H}
\end{align*}
is called the \emph{analysis operator} of the frame. Adjoint of $T_F^*$ is called the \emph{synthesis operator} of the $(\Omega,\mu)$-frame $\{F(w)\}_{w\in \Omega}$ and the operator $T_F^*T_F:\mathcal{H}\rightarrow \mathcal{H}$ is called \textit{frame operator} of $(\Omega,\mu)$-frame .

 \section{Dual of Retro  $(\Omega,\mu)$-Frames}

\begin{defn}
A system $\mathcal{F} \equiv (\{F(w)_{w\in \Omega}\},\Theta)$
$(\{F(w)_{w\in \Omega}\}\subset \mathcal{H },\Theta:\mathcal{Z}_d \rightarrow
\mathcal{H})$ is called a \emph{retro $(\Omega,\mu)$-frame} for
$\mathcal{H}$ with respect to an associated sequence space
$\mathcal{Z}_d$ if,
\begin{enumerate}[(i)]
\item  $\{\langle x,F(w)\rangle\}\in \mathcal{Z}_d$, for each $x\in \mathcal{H}$.
\vskip 2pt
\item There exist positive constants $(0< \mathrm{A}_0 \le \mathrm{B_0}< \infty)$ such that
$$\mathrm{A_0}\|x\|^2 \leq \int_\Omega|\langle x,F(w)\rangle|^2d\mu(w) \leq\mathrm{B_0}\|x\|^2, \ \text{for  each} \ x \in\mathcal{H}.$$
\vskip 2pt
\item $\Theta$ is a bounded linear operator such that
$\Theta(\{\langle x,F(w)\rangle\} = x, \ x \in \mathcal{H}$
\end{enumerate}

\end{defn}
The positive constant $\mathrm{A}_0,  \mathrm{B_0}$ are called the
\emph{lower} and \emph{upper} \emph{retro $(\Omega,\mu)$-frame  bounds} of
$(\{F(w)\}_{w\in \Omega},\Theta)$, respectively. As in case of
frames for Hilbert spaces, they are not unique.  The
operator $\Theta:\mathcal{Z}_d \rightarrow \mathcal{H}$ is called
\emph{retro pre $(\Omega,\mu)$-frame operator} (or simply $(\Omega,\mu)$-reconstruction operator)
associated with $\{F(w)\}_{w\in \Omega}$. A retro $(\Omega,\mu)$-frame
$\mathcal{F} \equiv (\{F(w)\}_{w\in \Omega},\Theta)$ is said to be \emph{exact}
if for each $w\in \Omega$ there exists no $(\Omega,\mu)$-reconstruction
operator $\Theta_w$ such that $(\{F(w)\}, \Theta_{w_0})_{w \ne w_0} $ is
retro $(\Omega,\mu)$-frame for $\mathcal{H}$.
\begin{lem}
Let $\mathcal{F}\equiv(\{F(w)\}_{w\in \Omega}, \Theta)$  be a retro $(\Omega,\mu)$-frame
for $\mathcal{H}$ . Then, $\mathcal{F}$  is exact if and only if
$F(w_0)\not\in [F(w)]_{w_0\neq w},$  for almost all $w\in \Omega$.
\end{lem}

\begin{lem}
Let $\mathcal{H}$ be a Hilbert space and let
$\left\{\{G(w)\}_{w\in \Omega}\right\}\subset \mathcal{H}$ be a sequence such
that $\left\{x \in \mathcal{H} :\langle x,G(w)\rangle= 0, \
\text{for almost all}\; w\in \Omega\right\} = \left\{0\right\}$. Then,
$\mathcal{H}$ is linearly isometric to the Hilbert  space
$\mathcal{Z} =\{\{\langle x,F(w)\rangle\}:x\in \mathcal{H} \}$, where the norm  is  given by

\begin{align*}
\int_\Omega|\langle x,G(w)\rangle|^2d\mu(w)=\|x\|^2
, \ x \in \mathcal{H}.
\end{align*}
\end{lem}

\begin{defn}
Let $\mathcal{G} \equiv (\{\langle x,G(w)\rangle\},\Theta)$ be a retro $(\Omega,\mu)$-frame
for $\mathcal{H}$ and let  $\Omega_0$  be a fixed subset of
$\Omega$. A system $\{\langle x,G(w)\rangle\} \subset \mathcal{H}$ is
called \emph{dual retro $(\Omega,\mu)$-frame} of $\mathcal{F}$ (or simply
$\Omega_0$-dual of $\mathcal{F}$) if
\begin{enumerate}[(i)]
\item $\langle G(w_0),F(w)\rangle= \delta_{w_0, w}$, for almost all $ w \in \Omega$ and  for almost all $w_0 \in
\Omega \setminus \Omega_0$,

\item there exists a reconstruction operator $\theta^\dag$ such that
$\mathcal{G}\equiv (\{G(w)\}, \Theta^\dag)$ is a retro $(\Omega,\mu)$-
frame for $\mathcal{H}.$
\end{enumerate}
\end{defn}

If $\mathcal{G}\equiv (\{G(w)\}, \Theta^\dag)$ is a dual of
$\mathcal{F} \equiv (\{F(w)\},\Theta)$, then we say that
$\mathcal{F}$ admits a dual with respect to the system
$\{G(w)\}$. In short we say that $(\mathcal{F}, \mathcal{G})$ is
a \emph{retro dual pair}.

\begin{rem}
In Definition $3.1$, if  $\Omega_0 = \emptyset$, then $\mathcal{G}$
is called the strong dual of $\mathcal{F}$.
\end{rem}

\begin{exa}
Consider the measure space $\mathcal{H}=L^2(\Omega)$ with counting measure, where $\Omega= \mathbb{N}$
 Let $\{F(w)\} = \{\chi_w\}$ be an
orthogonal basis for a Hilbert space $\mathcal{H}$.\\
Then
\begin{align*}
\int_\Omega|\langle x,F(w)\rangle|^2d\mu(w)=\|x\|^2,\ x
\in \mathcal{H}.
\end{align*}
Define $\Theta: \ell^2 \rightarrow \mathcal{H}$ by
$\Theta(\{\langle x,F(w)\rangle\}) = x $. Then, $\Theta$ is
reconstruction operator  such that
$\mathcal{F}\equiv(\{F(w)\},\Theta)$ is a normalized tight exact
retro $(\Omega,\mu)$-frame for $\mathcal{H}$ with respect to
$\mathcal{Z}_d = \ell^2$.
\end{exa}
 Let $\mathcal{F}$ be an exact retro $(\Omega,\mu)$-frame for
$\mathcal{H}$. Then, in general, $\mathcal{F}$ has no retro dual
$(\Omega,\mu)$-frame. The following example provides the existence of an exact
retro $(\Omega,\mu)$-frame $\mathcal{F}$ which does not admits a dual retro
$(\Omega,\mu)$-frame.

\begin{exa}
Consider the measure space $\mathcal{H}=L^2(\Omega)$ with counting
measure, where $\Omega= \mathbb{N}$.  Let $\{\chi_w\}$ be an
orthonormal basis of $\mathcal{H}$.\\ Define $\{F(w)\} \subset
\mathcal{H}$ by
\begin{align*}
 F(w) = \chi_{w+1}+\chi_1,\ w
\in \Omega.
\end{align*}
Let  $\mathcal{Z}= \{{\{\langle x,F(w)\rangle\}: x \in
\mathcal{H}}\}$. Then, by using Lemma $3.3$, $\mathcal{Z}$ is a
Hilbert space with norm given by
\begin{align*}
\int_\Omega |\langle x,F(w)\rangle|^2d\mu(w)=\|x\|^2,\ x \in
\mathcal{H}.
\end{align*}
 Define $\Theta: \mathcal{Z} \rightarrow \mathcal{H}$ by
\begin{align*}
\Theta(\{\langle x,F(w)\rangle\})= x, \ x \in \mathcal{H}.
\end{align*}
Then, $\Theta$ is a bounded linear operator such that
$\mathcal{F}\equiv(\{F(w)\},\Theta)$ is a retro $(\Omega,\mu)$-frame for
$\mathcal{H}$ with respect to $\mathcal{Z}$.

To show $\mathcal{F}$ is exact. Choose $G(w)= \chi_{w+1}, w
\in \Omega$. Then, $\{G(w)\} \subset \mathcal{H}$ and we
observe that
\begin{align*}
\langle G(w_0),F(w)\rangle= \delta_{w_0,w},\ \text{for almost all}\ w_0 , w \in \Omega.
\end{align*}
 Therefore,
$G(w_0)\not\in [F(w)]_{w\neq w_0},$ for almost all $w_0\in \Omega$. Thus,
by using Lemma $3.2$, we conclude that $\mathcal{F}$ is an exact
retro $(\Omega,\mu)$-frame for $\mathcal{H}$. Note the if $\mathcal{F}$ is
exact, then the condition $(i)$ given in definition $3.1$ is
satisfied for $\Omega_0 = \emptyset$.

To show that the condition $(ii)$ given in Definition $3.1$ is not
satisfied. Let  $\Theta^\dag$ be the reconstruction operator such
that $(\{F(w)\}, \Theta^\dag)$ is a retro $(\Omega,\mu)$-frame for
$\mathcal{H}$. Let $A^0, \ B^0$ be a choice of  retro bounds
for $(\{G(w)\}, \Theta^\dag)$.\\ Then
\begin{align}
A^0\|y\|^2 \leq \int_\Omega|\langle y,G(w)\rangle|^2d\mu(w)
\leq B^0 \|y\|^2, \ \text{for all}\ y
\in\mathcal{H}.
\end{align}
In particular for $y= \chi_1,$ we have $\langle y,G(w)\rangle
=0, $ for all $w \in \Omega$. Therefore, by using retro $(\Omega,\mu)$-frame
inequality  $(3.1)$, we obtain $y=0,$ a contradiction. Hence
$\mathcal{F}$ has no dual retro $(\Omega,\mu)$-frame.
\end{exa}
\begin{rem}
We observe that may have a retro $(\Omega,\mu)$-frame does not admit strong  duality, but
it may have a dual pair.
\end{rem}

The following example defend Remark $3.8$
\begin{exa}
Let $\mathcal{F}\equiv(\{F(w)\},\Theta)$ be a retro $(\Omega,\mu)$-frame
for $\mathcal{H}$ given in Example $3.7$. Choose $G(w)=
\chi_w, w \in \Omega$ and let $\mathcal{Z}^{0}= \{
\{\langle y,G(w)\rangle\}: y \in \mathcal{H}\}$. Then,
$\mathcal{Z}^{0}$ is a Hilbert space of sequences of scalars with
norm given by
 \begin{align*}
\int_\Omega|\langle y,G(w)\rangle|^2d\mu(w)= \| y\|^2, \   y \in \mathcal{H} .
\end{align*}
 Therefore, $\Theta^\dag:
\{\langle y,G(w)\rangle\}\rightarrow y $ is a bounded linear
operator from $\mathcal{Z}^{0}$ onto $\mathcal{H}$ such that
$\mathcal{G}_0 \equiv (\{G(w)\},\Theta^\dag)$ is a retro $(\Omega,\mu)$-frame for $\mathcal{H}$ with respect to
$\mathcal{Z}^{0}$, with bounds $A = B = 1$. Choose $\Omega_0 = \{w\}$.
Then, $\langle G(w_0),F(w)\rangle= \delta_{w_0, w}$, for almost all $ w \in \Omega$
and for almost all $w_0 \in \Omega \setminus \Omega_0$. Therefore,
$(\mathcal{F}, \mathcal{G}_0)$ is  a retro dual pair, which is not
strong.
\end{exa}

\begin{rem}
By Example $3.9$ we observe that  a retro Banach frame $\mathcal{F}$
for $\mathcal{H}$ may have dual pair but does not admit a strong
dual (even) with respect to an orthonormal basis for
$\mathcal{H}$.

\end{rem}

The following proposition provides  sufficient conditions for the
existence of  the  dual retro $(\Omega,\mu)$-frames.  It is sufficient to
prove the result for the existence of strong dual retro $(\Omega,\mu)$-frame. We can extend the same construction for the existence of
arbitrary dual pair.
\begin{prop}
Let $\mathcal{F} \equiv (\{F(w)\},\Theta)$ be a retro $(\Omega,\mu)$-frame
for $\mathcal{H}$. Then, $\mathcal{F}$  admits a strong dual retro
$(\Omega,\mu)$-frame if there exists a system $(\{G(w)\} \subset
\mathcal{H}$ such that $\langle G(w_0),F(w)\rangle= \delta_{w_0,w},$ for almost all
$w_0 , w \in \Omega$ and there exists an injective closed linear
operator
 $y \rightarrow \{\langle y,G(w)\rangle\}$ with closed
range from $\mathcal{H}$ to $\mathcal{\widehat{Z}}_0$, where
$\mathcal{\widehat{Z}}_0$ is some Hilbert space of scalar valued
sequences. These conditions are not necessary.
\end{prop}
\proof Let  $U:\mathcal{H} \rightarrow \mathcal{\widehat{Z}}_0$
be given by $U(y)= \{\langle y,G(w)\rangle\}, y \in
\mathcal{H}$. Then, by hypothesis, $U$ is injective and closed
linear operator  with closed range $R(U)$. Therefore,  $U^{-1}: R(U)
\subset \mathcal{\widehat{Z}}_0 \rightarrow \mathcal{H}$ is
closed [13]. Thus, by using Closed Graph Theorem [13], there exists
a constant $c>0$ such that
\begin{align*}
\|U(y)\|\geq c \|y\|,\ \text{for all} \ y
\in \mathcal{H}.
\end{align*}
 Let, if possible, there exists no reconstruction
operator $\Theta^\dag$ such that $(\{G(w)\},\Theta^\dag)$ is
a retro $(\Omega,\mu)$-frame for $\mathcal{H}$. Then, by using
Hahn-Banach Theorem there exists a non-zero functional $y_0
\in \mathcal{H}$ such that $\langle y_0,G(w)\rangle=0,$ for all
$w \in \Omega.$\\ Therefore
\begin{align*}
 0= \|U(y_0)\|\geq c \|y_0\|.
\end{align*}
 This
gives $y_0=0,$ a contradiction. Therefore, there exists a
reconstruction operator $\Theta^\dag$ such that
$(\{G(w)\},\Theta^\dag)$ is a retro $(\Omega,\mu)$-frame for
$\mathcal{H}$ with respect to some associated Banach space of
scalar valued sequences. Hence $\mathcal{F}$ admits a dual retro $(\Omega,\mu)$-frame for the underlying space.

To show that the conditions  are not necessary. Let $\mathcal{H}$ be
the measure space given in Example $3.7$. Let $F(w) = w^2\chi_w, w
\in \Omega$, where $\{\chi_w\}$ is an orthonormal basis for
$\mathcal{H}$. Then, there exists a bounded linear operator $\Theta$
such that $\mathcal{F} =(\{F(w)\}, \Theta)$ is a retro $(\Omega,\mu)$-frame for $\mathcal{H}$.\\
Choose $G(w)= \frac{1}{w^2}\chi_w, w\in \Omega$. Then,
$\langle G(w_0),F(w)\rangle= \delta_{w_0,w}$ for almost all $w_0, w \in \Omega$. By
the nature of the system $\{G(w)\}$, we conclude that there
exists a reconstruction operator $\Theta^\dag$ such that
$\mathcal{G}\equiv (\{G(w)\}, \Theta^\dag)$ is a retro $(\Omega,\mu)$-frame for $\mathcal{H}$ with respect to
$\mathcal{\widehat{Z}}_{o} = \ell^2$. Thus, $(\mathcal{F},
\mathcal{G})$ is a dual pair.

  Define $\widetilde{\Theta}:
\mathcal{H}\rightarrow \mathcal{\widehat{Z}}_{o} (= \ell^2)$ by
\begin{align*}
\widetilde{\Theta}(y)= \{\langle y,G(w)\rangle\}=
\{\frac{\xi_w}{w^2}\}, \ y = \{\xi_w\} \in \mathcal{H}.
\end{align*}
 It can be verified that  the range of the operator
$\widetilde{\Theta}$ is not closed.
\endproof

The following theorem gives  necessary and sufficient condition for
a given retro $(\Omega,\mu)$-frame to admit its dual.
\begin{thm}
Let $\mathcal{F} \equiv(\{F(w)\}, \Theta)$ be a retro $(\Omega,\mu)$-frame
 for $\mathcal{H}$. Then, $\mathcal{F}$ has
a dual retro $(\Omega,\mu)$-frame if and only if there exists a system
$\{G(w)\} \subset \mathcal{H}$  such that $\langle G(w_0),F(w)\rangle=
\delta_{w_0, w}$, for almost all $ w \in \Omega$, $w_0 \in
\Omega \setminus \Omega_0$ and dist$( x, L_w)\rightarrow 0$
as $w\rightarrow\infty$, for all $ x\in\ \mathcal{H}$, where
$L_w= [G(w_1),G(w_2),....,G(w)]$, for all $w \in
\Omega$.

\end{thm}
\proof Suppose first that $\mathcal{F} \equiv(\{F(w)\}, \Theta)$
has  a dual retro $(\Omega,\mu)$-frame. Then, by definition we can find a
system $\{G(w)\} \subset \mathcal{H}$  such that
$\langle G(w_0),F(w)\rangle= \delta_{w_0, w}$, for almost all $ w \in \Omega$, $w_0 \in \Omega \setminus \Omega_0$ and  a reconstruction
operator  $\Theta^\dag$ such that $\mathcal{G}\equiv
(\{G(w)\}, \Theta^\dag)$ is a retro $(\Omega,\mu)$-frame for
$\mathcal{H}.$ Let $A_0$ and $ B_0$ be a choice of bounds for
$\mathcal{G}$.\\
Then
\begin{align}
A_0\|y\|^2\leq\int_\Omega|\langle y,G(w)\rangle|^2d\mu(w)\leq
B_0\|y\|^2,\quad \text{for all } \
y\in \mathcal{H}.
\end{align}
Suppose that the condition  $\text{dist}( x, L_w)\rightarrow0$ as
$w\rightarrow\infty$, for all $ x \in\ \mathcal{H}$, is  not
satisfied. Then, there exists a non zero functional $y_0 \in
\mathcal{H}$ such that
\begin{align*}
\lim_{w\rightarrow \infty} \text{dist}(y_0, L_w) \ne 0.
\end{align*}
 Note that dist$(y_0, L_w)$ $\geq$
dist$(y_0, L_{w+1})$ for all $w \in \Omega$. This is
because $\{L_w\}$ is a nested system of subspaces and by definition
of distance of a point from a set. Now $\{\text{dist}(y_0,
L_w)\}$ is bounded below monotone decreasing sequence. So,
$\{\text{dist}(y_0, L_w)\}$ is convergent and its limit is a
positive real number, since otherwise dist$( y_0,
L_w)\rightarrow 0$ as $w\rightarrow\infty$ (which is not possible).
Let $\lim_{w\rightarrow \infty} \text{dist}( y_0, L_w) = \xi >
0.$

Choose $D= \bigcup_{w}L_w$. Then, by using the fact that
$\text{dist}(y_0, D) = \inf\{\text{dist}(y_0, L_w)\}$, we
obtain
\begin{align}
\text{dist}(y_0, D) \geq \xi > 0.
\end{align}

 Now we show that $y_0 \not\in
\overline{D}$. Let if possible, $y_0 \in \overline{D}$. Then,
we can find a sequence $\{\zeta_w\} \subset D$ such that
$\text{dist}(\zeta_w, y_0) \rightarrow 0 \ \text{as}\ w
\rightarrow \infty.$

By using $(3.3)$, we have dist$(y_0, D) \geq \xi$. Therefore,
dist$(\zeta_w, y_0) \geq \xi
>0.$ This is a contradiction to the fact that dist$(\zeta_w,
y_0) \rightarrow 0$ as $w \rightarrow \infty$. Hence $y_0
\not\in \overline{D}$. Thus, by using Hahn-Banach theorem, there
exists a non zero functional $z_0 \in \mathcal{H}$ such
that
\begin{align*}
 \langle z_0,G(w)\rangle= 0,\ \text{for all} \ w\in \Omega.
\end{align*}
  Therefore,
 by using retro $(\Omega,\mu)$-frame inequality $(3.2)$, we have $z_0=0$, a contradiction. Hence
$\text{dist}( x, L_w)\rightarrow0$ as $w\rightarrow\infty$, for
all $ x\in\ \mathcal{H}$.

To prove the converse part, assume that there exists a system
$\{G(w)\} \subset \mathcal{H}$ is such that
\begin{align*}
\langle G(w_0),F(w)\rangle= \delta_{w_0, w},\ \text{for almost all}\  w \in \Omega
\ \text{and for almost all} \ w_0 \in \Omega \setminus \Omega_0,
 \intertext{and}
\text{dist}( x, L_w)\rightarrow 0\ \text{as}\
w\rightarrow\infty,\ \text{for all} \  x \in \mathcal{H}.
\end{align*}
 Then, in particular, for each $\epsilon > 0$ and
for each $x \in \mathcal{H}$, we can find a $G(w_0)$ from
some $L_w$ such that
\begin{align*}
\|x- G(w)\| < \epsilon.
\end{align*}
Therefore, by using Lemma $3.3$,  $\mathcal{Z}= \{\{\langle y,G(w)\rangle\}: y \in \mathcal{H}\}$ is a
Hilbert space of sequences of scalars with norm given by
 \begin{align*}
\int_\Omega|\langle y,G(w)\rangle|^2d\mu(w)= \|y\|^2, \  y \in\mathcal{H} .
\end{align*}
Define $\Theta_0: \mathcal{Z}\rightarrow \mathcal{H}$ by
 \begin{align*}
\Theta_0(\{\langle y,G(w)\rangle\}) = y, \   y
\in\mathcal{H}.
\end{align*}
Then, $\Theta_0$ is a bounded linear operator such that
$(\{G(w)\}, \Theta_0)$ is a retro $(\Omega,\mu)$-frame for
$\mathcal{H}$ with respect to $\mathcal{Z}$. Hence
$\mathcal{F}$ has a dual retro $(\Omega,\mu)$-frame.
\endproof

The existence of duality with respect to certain sequence space is
also interested. The following theorem provides the necessary and
sufficient condition for an exact retro $(\Omega,\mu)$-frame to admit a dual
frame with respect to a given sequence space.

\begin{thm}
Let $\mathcal{F}  \equiv (\{F(w)\}, \Theta)$ be a retro $(\Omega,\mu)$-frame for $\mathcal{H}$ with respect to $\mathcal{Z}_\mathfrak{d}$
and let $\mathcal{A}_d = \{\{\langle y,G(w)\rangle\}: y \in
\mathcal{H} \}$. Then, $\mathcal{F}$ has a dual retro $(\Omega,\mu)$-frame with respect to the sequence space $\mathcal{A}_d$ if and only
if there exists a system $\{G(w)\} \subset \mathcal{H}$  such
that $\langle G(w_0),F(w)\rangle= \delta_{w_0, w}$, for almost all $ w \in \Omega$, $w_0 \in \Omega \setminus \Omega_0$ and the analysis
operator $\mathcal{U}: y \rightarrow
\{\langle y,G(w)\rangle\}$ is a bounded below continuous linear
operator from $\mathcal{H}$
 onto $\mathcal{A}_d$.
\end{thm}
\proof Suppose first that $\mathcal{F}$ has a dual retro $(\Omega,\mu)$-frame $\mathcal{G}$ with respect  to $\mathcal{A}_d$. Then, there
exists a reconstruction operator $\Theta^\dag$ such that
$\mathcal{G}\equiv (\{G(w)\}, \Theta^\dag)$ is a retro $(\Omega,\mu)$-frame for $\mathcal{H}$ with respect to $\mathcal{A}_d$.
Therefore, there are positive constants $\mathrm{a_0}$ and $
\mathrm{b_0}$ such that
\begin{align}
\mathrm{a_0}\|y\|^2 \leq
\int_\Omega|\langle y,G(w)\rangle|^2d\mu(w)
\leq\mathrm{b_0}\|y\|^2, \ \text{for \ each} \ y
\in\mathcal{H}.
\end{align}
Now consider the analysis operator $\mathcal{U}: \mathcal{H}
\rightarrow \mathcal{A}_d$ which is given by
\begin{align*}
\mathcal{U}(y) = \langle y,G(w)\rangle, \ y
\in\mathcal{H}.
\end{align*}
Then, linearity and ontoness of $\mathcal{U}$ is obvious. By using
upper retro $(\Omega,\mu)$-frame inequality in $(3.4)$, we have
\begin{align*}
\int_\Omega|\mathcal{U}(y)|^2d\mu(w)
\leq\mathrm{b_0}\|y\|^2, \ \text{for \ each} \ y
\in\mathcal{H}.
\end{align*}
Therefore, $\|\mathcal{U} \| \leq \mathrm{b_0}.$ Hence $\mathcal{U}$
is continuous. Similarly, by using lower retro frame inequality in
$(3.4)$, we have
\begin{align*}
\int_\Omega|\mathcal{U}(y)|^2d\mu(w)\geq \mathrm{a_0}
\|y\|^2.
\end{align*}
Hence $\mathcal{U}$ is bounded below.\\
For the reverse part, assume that  $\mathcal{U}$ is bounded below.
Then, using Lemma $3.3$, $\mathcal{A}_d$ is a Hilbert space with the
norm given by
\begin{align*}
\int_\Omega|\langle y,G(w)\rangle|^2d\mu(w) =
\|y\|^2, \  y \in \mathcal{H}.
\end{align*}
 Define $\Theta^\dag:\mathcal{A}_d\rightarrow \mathcal{H} $ by
$\Theta^\dag(\langle y,G(w)\rangle) = y$.
 Then, $\Theta^\dag$ is a bounded linear operator such that $(\{G(w)\}, \Theta^\dag)$ is a retro $(\Omega,\mu)$-frame
for $\mathcal{H}$ with respect to $\mathcal{A}_d$. The theorem
is proved.
\endproof

\mbox{}

\begin{thebibliography}{99}\baselineskip11pt

\bibitem{AAG}
S. Ali, Tawreque, J.P. Antonie, J.P. Gazeau, continuous frames in Hilbert spaces, \emph{Ann. Physics},222(1) (1993), pp.38--88

\bibitem{CK}
 P.G.Casazza and G. Kutynoik,
Finite Frames, Birkh$\ddot{a}$user, 2012.

\bibitem{PC} P.G. Casazza,
The art of frame theory, \emph{Tawanese J. Math.}, 4(2) (2000),
129--201




\bibitem{CKL}
P.G.Casazza, G.Kutyniok and M.C.Lammers, Duality principles in frame
thoery, \emph{J. Fourier Anal. Appl.}, 10 (2004), 383--408.


\bibitem{C}
O. Christensen, Frames and bases (An introductory course),
Birkh\"auser, Boston (2008).


\bibitem{DGM} I.Daubechies, A. Grossmann and Y. Meyer,
Painless non-orthogonal expansions, \emph{J. Math. Phys.} 27 (1986),
1271--1283.

\bibitem{DS}  R.J. Duffin and A.C. Schaeffer,
A class of non-harmonic Fourier series, \emph{Trans. Amer. Math.
Soc.}, 72 (1952), 341--366.


\bibitem{FG}  H.G. Feichtinger and K. Gr\"ochenig,
A unified approach to atomic decompositons via inegrable group
representations, \emph{Lecture Notes in Mathematics}, 1302
(Springer, Berlin, 1988), 52--73.


\bibitem{GH}
J. P. Gabardo, D. Han, Frames associated with measurable spaces, \emph{Adv. Comput. Math.,} 18 (2003), pp. 127--47.


\bibitem{G} K. Gr\"ochenig,
Describing functions: Atomic decompositions versus frames,
\emph{Monatsh. Math.}, 112, (1991), 1--41.

\bibitem{HL}  D. Han and D.R. Larson,
Frames, bases and group representations, \emph{Mem. Amer. Math.
Soc.}, 147 (697) (2000), 1--91.


\bibitem{H} C. Heil,
A basis theory primer, Birkh\"auser (expanded edition)(1998).


\bibitem{HH}
H.Heuser, Functional Analysis, John Wiley and Sons, New York (1982).

\bibitem{J} P.K. Jain, S.K. Kaushik and L.K. Vashisht,
Banach frames for conjugate Banach spaces,
\emph{Zeitschrift f\"ur Analysis und ihre Anwendungen}, \textbf{23} (4) (2004), 713--720.

\bibitem{JK} P.K. Jain, S.K. Kaushik and L.K. Vashisht,
On perturbations of Banach frames,
\emph{International Jour. of Wavelet, Multiresolution and Information Processing} (IJWMIP), \textbf{4} (3) (2006), 559--565.



\bibitem{K}
G. Kaiser, A Friendly Guide to Wavelets, Birkh\"{a}user, Boston (1994), 300p.


\bibitem{KV}
 S.K. Kaushik, L.K. Vashisht, S.K. Sharma, Some results concerning frames associated with measurable spaces, \emph{TWMS J. Pure Appl. Math.,}
4 (1) (2013), pp. 52--60.


%\bibitem{18} L.K. Vashisht, On retro Banach frames of type $P$, \emph{Azerb. J.
%Math.}, 2 (1) (2012), 82--89.
%
%
%\bibitem{19} L.K. Vashisht, On $\Phi$-Schauder frames, \emph{TWMS J.  App. and Eng.
%Math.(JAEM)}, 2 (1) (2012), 116-120.

\bibitem{Y}
R. Young, On complete biorthogonal systems, \emph{Proc. Amer. Math.
Soc.}, 83 (3) (1981), 537-540.


\bibitem{Y1}
R. Young,
A introduction to non-harmonic Fourier series, Academic Press, New
York  (revised first edition 2001).


\end{thebibliography}
\end{document}